\ifx\shlhetal\undefinedcontrolsequence\let\shlhetal\relax\fi

\documentclass[11pt]{amsart}

\usepackage{amsmath}
\usepackage{amssymb}

\newtheorem{theorem}{Theorem}[section]
\newtheorem{claim}[theorem]{Claim}

\newtheorem{lemma}[theorem]{Lemma}

\newtheorem{corollary}[theorem]{Corollary}

\newtheorem{conjecture}[theorem]{Conjecture}

\theoremstyle{definition}
\newtheorem{definition}[theorem]{Definition}

\theoremstyle{remark}
\newtheorem{remark}[theorem]{Remark}

\newcount\skewfactor
\def\mathunderaccent#1#2 {\let\theaccent#1\skewfactor#2
\mathpalette\putaccentunder}
\def\putaccentunder#1#2{\oalign{$#1#2$\crcr\hidewidth
\vbox to.2ex{\hbox{$#1\skew\skewfactor\theaccent{}$}\vss}\hidewidth}}

\def\smallbox#1{\leavevmode\thinspace\hbox{\vrule\vtop{\vbox
   {\hrule\kern1pt\hbox{\vphantom{\tt/}\thinspace{\tt#1}\thinspace}}
   \kern1pt\hrule}\vrule}\thinspace}

\newcommand{\Then}{{\underline{Then}}}

\def\qedref#1{$\qed_{\reforiginal{#1}}$}

\title{simple wedge points}
\author{Shimon Garti}
\address{Institute of Mathematics
 The Hebrew University of Jerusalem
 Jerusalem 91904, Israel}
\email{shimon.garty@mail.huji.ac.il}


\begin{document}
\let\labeloriginal\label
\let\reforiginal\ref

\begin{abstract}
Let $V$ be a finite set of points in the plane, not contained in a line. Assume $|V| = n$ is an odd number, and $|L \cap V| \leq 3$ for every line $L$ which is spanned by $V$. We prove that every simple line $L_{a,b}$ in $V$ creates a simple wedge (i.e., a triple $\{a,b,c\} \subseteq V$ such that $L_{a,b}$ and $L_{a,c}$ are simple lines). We also show that both restrictions on $V$ (namely $|V|$ is odd and $|L \cap V| \leq 3$) are needed.

\par \noindent We conjecture, further, that if $|V| = n$ is an odd number then $V$ contains a simple wedge, even if $V$ is not $3$-bounded. We introduce a method for proving this, which gives (in this paper) partial results.
\end{abstract}

\maketitle

\newpage

\section{introduction}

Let $V$ be a finite set of points in the plane. Trying to analize the lines which are spanned by the points of $V$, we may assume that the size of $V$ is at least $3$, and $V$ is not contained in a line. We call a set of points with these properties \emph{an interesting set}.

The Gallai-Sylvester Theorem asserts that every interesting set $V$ spans a line which contains exactly two of the points in $V$. We call this line \emph{a simple line}, and if the points are $a,b$ we denote it by $L_{a,b}$. So the Gallai-Sylvester Theorem shows that there is at least one simple line. An elegant proof of this theorem (due to Kelly) can be found in \cite{MR1801937}.

In most cases, there is more than one simple line. The amount of simple lines was investigated in the paper of Dirac, \cite{MR0043485}, and he conjuctured that this number is at least $\frac{1}{2}n$. In fact, we know that in many cases $\frac{1}{2}n$ is also an upper bound, for instance when $n$ is even, see \cite{MR1194036}. Despite some efforts (and partial results, see \cite{MR0097014} and \cite{MR1194036}), Dirac's conjecture is still open.

Our main concern is the existence of a structure called \emph{a simple wedge}. $\{ a, b, c \} \subseteq V$ is a simple wedge if $a$ is the intersection of two simple lines spanned by $V$. It means that there are two points $b,c \in V$ such that both $L_{a,b}$ and $L_{a,c}$ are simple. If $n \in \{3,4,5 \}$ then any configuration of $n$ points (not contained in a line) has a simple wedge. But this does not hold for every $n$, and one can construct easily a configuration of $6$ points in the plane, not contained in a line, with no simple wedge (we shall describe this construction at the end of the first section).

So one problem is whether an interesting set of points contains a simple wedge. We can raise a stronger demand, that \emph{from every} simple line $L_{a,b}$ it will be possible to create a simple wedge. As we shall see, if $V$ is assumed to be $3$-bounded (i.e., $|L \cap V| \leq 3$ for every $L$ spanned by $V$) then we can give a positive answer when the size of $V$ is an odd number.

One may wonder if the restrictions on $V$ are essential. We shall see that, in a meaning, both of them are needed. On one hand, we can describe a set of $6$ points with no simple wedge at all, despite the fact that this set is $3$-bounded. Such an example is also available for every even number $n>4$, if $n = 2 {\rm mod} 4$. It means that the restriction on the odd size of $V$ is vital (but might be improvable). On the other hand, we can build a set of $9$ points, not $3$-bounded, with a simple line $L_{a,b}$ in it that does not create a simple wedge. It follows that the demand of $3$-boundedness is also essential.

Still, we can remove the strong requirement of creating a simple wedge out of every simple line, and ask what are the conditions on the size of $V$ so that every set of this size contains a simple wedge. In particular, we may ask if this set is infinite. The following is plausible:

\begin{conjecture}
\label{mt}
Let $V$ be a set of $n$ points in the plane, $n$ is an odd number, and $V$ is not contained in a line. \newline 
\Then\ $V$ spans a simple wedge.
\end{conjecture}

We shall introduce a partial result concerning this conjecture. Clearly, it is connected to the problem of finding the number of simple lines. Observe that if $|V|=n$ and $V$ spans more than $\frac{n}{2}$ simple lines, then there exists a simple wedge in $V$. But as indicated above, since $\frac{n}{2}$ is pretty closed to the upper bound (in many cases, it can be proved that the upper bound does not exceed $\frac{n}{2}$), we cannot lean just on the number of simple lines for our conjecture.

In the light of having $\frac{n}{2}$ (or less) as an upper bound, conjecture \ref{mt} becomes much more interesting. Despite the fact that the pure combinatorial argument (based on the amount of simple lines) does not force the existence of a simple wedge, any configuration of $n$ points (when $n$ is odd!) must contain one. If this conjecture is true, it means that the existence (and the number) of simple lines (and more complicated structures) depends on the geometry of $\mathbb{R}^2$. It might explain some of the difficulty to generalize the Gallai-Sylvester theorem to other settings (see \cite{MR2060634} for metric spaces).

For proving conjecture \ref{mt} we suggest the following scheme. Given a set $V$ which is not $3$-bounded, we try to use a slight perturbation on the points of $V$. Our purpose is to produce $V'$ which is $3$-bounded, and then to find a simple wedge in $V'$ (by virtue of the main result in section one here). We must take care of the accurate perturbation, so that this simple wedge in $V'$ will be also a simple wedge in $V$.

To exploit this idea, we might have to employ a nested process of induction (on the degree of boundedness, the nunber of $k$-bounded lines, and so on). We shall prove, in this paper, just a basic step; but this step demonstrates the general method. We shall see that if $V$ is $3$-bounded except of one line $L$ such that $|L \cap V| = 4$ then $V$ contains a simple wedge. We emphasize that in this case $V$ may fail the strong theorem of creating a simple wedge out of every simple line. Anyway, we do not know how to prove the general case.

The paper is arranged in two sections. The first section includes the main theorem, applied to $3$-bounded sets of points. The second section is devoted to the perturbation method, exemplified here in the specific case of having one line with more than $3$ points from $V$.

We indicate that the theorems below invite further investigation in several directions. For example, one can ask how many simple wedges do we have (as a function of $n$). The private cases of $n \in \{3,4,5\}$ yield the conjecture that there are always at least $3$ simple wedges (see corollary \ref{swfromeveryline} below). We can also try to analize more complicated and more general structures of simple lines, and so on. We hope to shed light on these questions in a subsequent work.

\newpage

\section{$3$-bounded sets}

\begin{definition}
\label{sw}
A simple wedge. \newline
Suppose $V$ is a set of points in the plane, not contained in a line. \newline
A simple wedge in $V$ is a triple $\{a,b,c\} \subseteq V$ such that the lines $L_{a,b}$ and $L_{a,c}$ are simple.
We call $a$ a simple wedge point.
\end{definition}

Let us try to draw the picture:

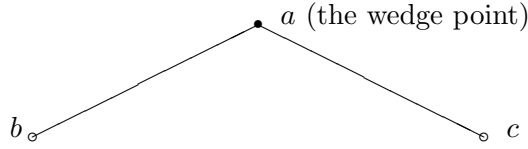
\setlength{\unitlength}{3cm}
\begin{picture}(1,1)
\put(1,0){\line(2,1){1}}
\put(1,0){\circle{0.04}}
\put(0.9,0){$b$}
\put(2,0.5){\circle*{0.04}}
\put(2.1,0.5){$a$ (the wedge point)}
\put(2,0.5){\line(2,-1){1}}
\put(3,0){\circle{0.04}}
\put(3.1,0){$c$}
\end{picture}

\medskip

So $a$ is a wedge point, and if the lines $L_{a,b}$ and $L_{a,c}$ are simple then $a$ is even a simple wedge point. Our main goal is to show that simple wedge points exist if $|V|$ is an odd number and the lines of $V$ contain at most $3$ points from $V$. We define:

\begin{definition}
\label{ellbounded}
$\ell$-boundedness. \newline 
Let $V$ be a set of points in the plane, $\ell$ a natural number. \newline 
$V$ is $\ell$-bounded if $|L \cap V| \leq \ell$ for every line $L$ which is spanned by the points of $V$.
\end{definition}

\begin{claim}
\label{mc}
The main claim. \newline 
Suppose $V$ is a finite set of points in the plane, not contained in a line, $n = |V|$ is an odd number and $V$ is $3$-bounded. \newline 
\Then\ $V$ spans a simple wedge.
\end{claim}

The main idea in proving this claim is using orbits of points, starting from a simple line of $V$. One can think of an orbit as a try to develop a simple wedge out of a simple line. The formal definition is as follows:

\begin{definition}
\label{orrbitss}
Orbits. \newline 
Let $V$ be a set of points in the plane, $a, b \in V$, and $L_{a, b}$ a simple line in $V$.
\begin{enumerate}
\item [$(a)$] $x = \langle x_1, \ldots, x_t \rangle$ is an orbit of $L_{a, b}$ when:
\begin{enumerate}
\item [$(\aleph)$] $x_i \notin \{a, b\}$ for every $1 \leq i \leq t$
\item [$(\beth)$] $x_i \neq x_j$ for every $\{i,j\} \subseteq \{1, \ldots, t-1\}$
\item [$(\gimel)$] $x_{2m + 1} \in L_{a, x_{2m}}$ for every relevant $m$
\item [$(\daleth)$] $x_{2m} \in L_{b, x_{2m - 1}}$ for every relevant $m$
\end{enumerate}
\item [$(b)$] $x$ is an open orbit if $t = 1$ or $t > 1$ and $x_1 \neq x_t$
\item [$(c)$] $x$ is a closed orbit if $t > 1$ and $x_1 = x_t$
\end{enumerate}
\end{definition}

The important thing to notice is that in an open orbit all the points are distinct, and in a closed orbit almost all the points are distinct, but the last point $x_t$ coincides with the first point $x_1$. A fundamental feature of orbits, for our proof, is the length of them:

\begin{definition}
\label{oorbitlg}
Length of orbits. \newline 
Let $x = \langle x_1, \ldots, x_t \rangle$ be an orbit for $L_{a,b}$.
\begin{enumerate}
\item [$(a)$] The length of $x$ is $t$ if $x$ is an open orbit, and $t-1$ if $x$ is a closed orbit
\item [$(b)$] $x$ is a maximal orbit for $L_{a,b}$ if there is no $y \in V$ so that $\langle x_1, \ldots, x_t, y \rangle$ is an orbit for $L_{a,b}$
\end{enumerate}
\end{definition}

\begin{remark}
\label{mmaxxexist}
\begin{enumerate}
\item [$(a)$] If $x = \langle x_1, \ldots, x_t \rangle$ is a closed orbit then $t>2$
\item [$(b)$] For every simple line $L_{a,b}$ and every $x_1 \in V$ not on this line, there exists a maximal orbit such that $x_1$ is its first point.
\end{enumerate}
\end{remark}

We proceed now to draw a connection between maximal open orbits and the existence of simple wedges:

\begin{lemma}
\label{maxopenorbit}
The characterization lemma. \newline 
Suppose $V$ is a set of points in the plane, not contained in a line. Let $L_{a,b}$ be a simple line in $V$. \newline 
One can create a simple wedge out of $L_{a,b}$ iff there is a maximal open orbit for this line.
\end{lemma}

\par \noindent \emph{Proof}.
\par \indent $\Rightarrow$ Assume that $\{a,b,c\} \subseteq V$ is a simple wedge for $L_{a,b}$ (without loss of generality $b$ is the wedge point, so the lines $L_{a,b}$ and $L_{b,c}$ are simple). It means that there is no $x \in V \setminus \{b,c\}$ which lies on the line $L_{b,c}$. Consequently, $\langle c \rangle$ is a maximal open orbit for the line $L_{a,b}$, as required.

$\Leftarrow$ Assume that $x = \langle x_1, \ldots, x_t \rangle$ is a maximal open orbit for the simple line $L_{a,b}$. If $t$ is an odd number then $L_{b, x_t}$ is a simple line. If $t$ is even then $L_{a, x_t}$ is a simple line. Either way, the triple $\{a,b,x_t\}$ establishes a simple wedge, so we are done.

\hfill \qedref{maxopenorbit}

\medskip

In the light of the former lemma, we would like to show that in every set of points $V$ such that $|V|$ is an odd number, there exists a maximal open orbit. The key point here is the following fact, that the length of a closed orbit (in a $3$-bounded set of points) is an even number:

\begin{lemma}
\label{eeven}
The even length lemma. \newline 
Let $V$ be a finite $3$-bounded set of points in the plane, not contained in a line. Suppose $L_{a,b}$ is a simple line in $V$, and let $x = \langle x_1, \ldots, x_t \rangle$ be a closed orbit for $L_{a,b}$. \newline 
\Then\ the length of $x$ (i.e., $t-1$) is an even number.
\end{lemma}

\par \noindent \emph{Proof}. \newline 
By the definition of length (for closed orbits), $t>2$ (see remark \ref{mmaxxexist}(a) above). Assume toward contradiction that $t-1$ is an odd number. By property $(\gimel)$ in definition \ref{orrbitss} of orbits, $x_t \in L_{b,x_{t-1}}$. Since $x$ is closed, $x_t = x_1$, hence $x_1 \in L_{b,x_{t-1}}$. It follows that the line $L_{b,x_{t-1}}$ is the line $L_{b,x_1}$, and we know that $x_2 \in L_{b,x_1}$. Let us call this line $L$.

Recall that $t>2$ hence $t-1 \geq 2$ so $t-1 > 2$ (by the assumption toward contradiction), so $x_2 \neq x_{t-1}$. Consequently, $b, x_1, x_2, x_{t-1} \in L \cap V$, in contrary to the assumption of the $3$-boundedness.

\hfill \qedref{eeven}

\begin{remark}
\label{ooddss}
The proof that the length of a closed orbit is an odd number cannot be improved. For every even number $n>2$, one can introduce a closed orbit of length $n$ above some simple line $L_{a,b}$.
\end{remark}

Let $x = \langle x_1, \ldots, x_t \rangle, y = \langle y_1, \ldots, y_r \rangle$ be two orbits. We say that $x$ and $y$ are disjoint if the sets $\{x_1, \ldots, x_t\}$ and $\{y_1, \ldots, y_r\}$ are disjoint. The following lemma describes the relationship between two orbits of the same simple line in $V$:

\begin{lemma}
\label{sseparation}
The separation lemma. \newline 
Let $V$ be as above (i.e., a finite $3$-bounded set of points in the plane, not contained in a line). Suppose $L_{a,b}$ is a simple line in $V$, and let $x = \langle x_1, \ldots, x_t \rangle$ be a closed orbit for $L_{a,b}$.
Suppose there is a point $y_1 \in V \setminus \{x_1, \ldots, x_t,a,b\}$ and let $y$ be an orbit which starts from $y_1$. \newline 
\Then\ $x$ and $y$ are disjoint.
\end{lemma}

\par \noindent \emph{Proof}. \newline 
Assume toward contradiction that $x$ and $y$ are not disjoint, and let $j$ be the first natural number so that $y_j = x_i$ for some $1 \leq i \leq t$. By the choice of $y_1$ it follows that $j>1$, hence $y_{j-1}$ is well defined. Call $L$ to the line between $y_j$ and $y_{j-1}$. Since $y$ is an orbit, either $a$ or $b$ lies on $L$, so without loss of generality it is $b$.

Notice that $L$ is also the line through $b,x_i$. There is a point $x' \in x$ such that $b$ lies on the line through $x',x_i$ (in fact, $x'$ is either $x_{i-1}$ or $x_{i+1}$). Again, this line is $L$ (since both $b$ and $x_i$ belong to this line), and in particular $x' \in L$. $x' \neq b$ and $x' \neq x_i = y_j$. Moreover, $x' \neq y_{j-1}$ by the minimality of $j$. Consequently, $b,x',y_{j-1},y_j \in L \cap V$ and these points are distinct. It follows that $|L \cap V| \geq 4$, a contradiction.

\hfill \qedref{sseparation}

Having the lemmas above, we can prove now the main claim. Recall that for every $3$-bounded $V$ whose size is an odd number, we have to find a simple wedge spanned by $V$.

\medskip 

\par \noindent \emph{Proof of the main claim \ref{mc}}: \newline 
Begin with a simple line $L_{a,b}$, spanned by $V$. Such a line exists by virtue of the Gallai-Sylvester theorem. $|V| \geq 3$, so choose $x_1 \in V \setminus \{a,b\}$. Let $x = \langle x_1, \ldots, x_t \rangle$ be a maximal orbit for $L_{a,b}$. If $x$ is an open orbit, we are done due to lemma \ref{maxopenorbit}. 

If $x$ is closed then the length of $x$ is an even number (according to lemma \ref{eeven}). The amount of points in $x \cup \{a,b\}$ is still an even number, so there exists $y_1 \in V \setminus \{x_1, \ldots, x_t,a,b\}$. Let $y = \langle y_1, \ldots, y_r \rangle$ be a maximal orbit for $L_{a,b}$. By the separation lemma \ref{sseparation}, $y$ is disjoint to $x$. Again, if $y$ is an open orbit then we can create a simple wedge, and if not then we can choose a new starting point $z_1$ to begin another orbit.

After a finite number of steps we shall remain with a maximal open orbit (recall that if we stay with a single point $c$, then $\langle c \rangle$ is an open orbit). By the characterization lemma \ref{maxopenorbit} we have a simple wedge in $V$, so the proof is complete.

\hfill \qedref{mc}

As a matter of fact, the proof above shows a little bit more:

\begin{corollary}
\label{swfromeveryline}
Many simple wedges. \newline 
Assume $V$ is a finite $3$-bounded set of points in the plane, not containd in a line, $|V|$ is an odd number. \newline 
\Then\ a simple wedge can be created from every simple line of $V$. \newline 
In particular, such $V$ contains at least $3n / 13$ simple wedges.
\end{corollary}

\par \noindent \emph{Proof}. \newline 
The first assertion follows from the main claim, and the second is due to \cite{MR1194036}.

\hfill \qedref{swfromeveryline}

We indicate that this corollary need not be true when the set $V$ is not $3$-bounded. An example with $9$ points is easy to construct. We give two examples, one of $6$ points with no simple wedge, and one of $9$ points with a simple line that does not create a simple wedge.

Define $a = (-2,0)$ and $b = (2,0)$. So far, the line $L_{a,b}$ is simple, and it will remain simple although we add more points. Set $x_1=(-1,2), x_2=(1,2), x_3=(0,4)$. At last, let $y=(0,\frac{4}{3})$. This is the intersection point between $L_{a,x_2}$ and $L_{b,x_1}$. Now let $V$ be $\{a,b,x_1,x_2,x_3,y\}$.

$|V|=6$, and there are exactly $3$ simple lines in $V$, namely $L_{a,b}$, $L_{x_1,x_2}$ and $L_{x_3,y}$. Evidently, no simple wedge exists in $V$. Notice that $V$ is $3$-bounded, but $|V|$ is not an odd number. This phenomenon is not a peculiar feature of the number $6$, and it can be done also for other natural numbers. In particular, it is possible to produce a $3$-bounded set of points of size $n$, for every even number $n>6$, so that there is a simple line $L_{a,b}$ with no simple wedge. The construction follows from remark \ref{ooddss}, by creating a closed orbit of length $n-2$ over a simple line $L_{a,b}$ (but there are other simple wedges in these constructions, when $n>6$).

Now we add more three points to $V$, aiming to create a set which is not $3$-bounded. Let $g_3$ be $(0,2)$ and $g_1=(-\frac{2}{3}, \frac{8}{3}), g_2=(\frac{2}{3}, \frac{8}{3})$. The point $g_3$ is the intersection between the line through $a,g_2$ and the line through $b,g_1$. Set $W = V \cup \{g_1,g_2,g_3\}$. 

The line $L_{a,b}$ is simple in $W$, but no other simple line comes out of $a$ or $b$. Notice that it happens despite the fact that $|W|$ is an odd number (but $W$ is not $3$-bounded, of course). Again, this is not a singularity of the number $9$. One can add more triples of the same type as $g_1,g_2,g_3$ and create larger sets of points with no simple wedge out of $L_{a,b}$.
These examples justify the assumptions on $V$ in the main claim of the paper.

\newpage 

\section{The perturbation method}

We try to describe, in this section, a possible way toward proving Conjecture \ref{mt}. The theorem below is only a private case (the corollary from this theorem is slightly more general), and the main point is the method of proof. Still, we do not know if it can be applicable to the most general case.

\begin{theorem}
\label{pperturb}
A single perturbation theorem. \newline 
Suppose $V$ is an interesting set, $n = |V|$ is odd. \newline 
Assume there is a line $L$ so that $|L \cap V| = 4$, and $|L' \cap V| \leq 3$ for every other line spanned by $V$. \newline 
\Then\ $V$ contains a simple wedge.
\end{theorem}

\par \noindent \emph{Proof}. \newline 
Let $L_{a,b}$ be a simple line in $V$ (so $L \neq L_{a,b}$). Since $a, b \in L$ is impossible, we may assume without loss of generality that $a \notin L$. If there is $c \in L$ so that $L_{a,c}$ is a simple line, we are done. If not, choose any $c \in L$, and notice that $|L_{a,c} \cap V| = 3$ (by our assumptions). Let $d$ be the third point (from $V$) which lies on $L_{a,c}$.

Choose $c' \in L_{a,b} \setminus \{a,b\}$, and replace $c$ by $c'$, i.e., set $V' = V \cup \{c'\} \setminus \{c\}$. $V'$ is $3$-bounded, and $|V'| = n$. By \ref{mc}, every simple line in $V'$ creates a simple wedge. In particular, the line $L_{a,d}$ (which is simple in $V'$, after the removal of the point $c$) creates a simple wedge.

If $a$ is a simple wedge point (with respect to $L_{a,d}$), then there exists $e \in V' \setminus \{a,d\}$ so that $L_{a,e}$ is a simple line. Notice that $e \notin L_{a,b}$ (since $L_{a,b}$ is not simple in $V'$ after adding $c'$). Consequently, $L_{a,e}$ is simple also in $V$, hence the triple $\{a,d,e\}$ establishes a simple wedge in $V$.

Suppose that no such $e$ exists in $V'$. It follows that $d$ is the simple wedge point (with respect to the line $L_{a,d}$). Let $d' \in V'$ be so that $L_{d,d'}$ is a simple line in $V'$ (which differs from $L_{a,d}$). As before, $d' \notin L_{a,d}$ hence $L_{d,d'}$ is a simple line also in $V$. Moreover, we can show that there is another point $d'' \in V$ such that $L_{d,d''}$ is a simple line in $V$.

For this, pick any point $p \in L_{d,d'} \setminus \{d,d'\}$, and set $W = V \cup \{p\} \setminus \{c\}$. As in the case of $V'$, $W$ is $3$-bounded and $|W| = n$ is an odd number. $L_{a,d}$ is simple in $W$ (after the removal of the point $c$), and $d$ must be the simple wedge point (since we eliminated the case of $a$ is the simple wedge point).

Now $L_{d,d'}$ is not simple in $W$ (after adding the point $\{p\}$), hence there exists $d'' \notin L_{d,d'}$ so that $L_{d,d''}$ is simple in $W$. Clearly, $L_{d,d''}$ is simple also in $V$ (recall that $d'' \notin L_{a,d}$ by its simplicity in $W$). We conclude that the triple $\{d,d',d''\}$ is a simple wedge in $V$, and the proof is complete.

\hfill \qedref{pperturb}

It is tempting to try some induction process, by applying the perturbation method for the general induction step. But there is an obstacle here. The first basic step of \ref{pperturb} goes back to the case of $3$-bounded set (after changing the location of one point from $L$). In this case, we know that \emph{every} simple line creates a simple wedge, and that was crucial in the proof above.

In the general case, a perturbation of one point does not create a $3$-bounded set. So the induction hypothesis is weaker (namely, there exists a simple wedge). Nevertheless, we strongly believe that conjecture \ref{mt} is correct. One way to generalize \ref{pperturb} is by using a perturbation of many points simultaneously. Denote the number of simple lines in $V$ by $s_V$. The following can be proved:

\begin{corollary}
\label{ppperturb}
Suppose $V$ is an interesting set, $n = |V|$ is odd. \newline 
Assume there is a line $L$ so that $|L \cap V| = k$ when $4 \leq k \leq s_V+3$, and $|L' \cap V| \leq 3$ for every other line spanned by $V$. \newline 
\Then\ $V$ contains a simple wedge.
\end{corollary}

\par \noindent \emph{Proof}. \newline 
Let $\{L_1, \ldots, L_{s_V}\}$ be an enumeration of the simple lines in $V$, without repetitions. Let $\{v_1, \ldots, v_k\}$ be an enumeration of $L \cap V$, and without loss of generality $k = s_V+3$. For every $1 \leq j \leq s_V$ choose $c_j \in L_j$ which differs from the points of $L_j \cap V$. Set $V' = V \cup \{c_1, \ldots, c_{s_V}\} \setminus \{v_1, \ldots, v_{s_V}\}$. Now proceed as in the proof of theorem \ref{pperturb}.

\hfill \qedref{ppperturb}

\newpage

\bibliographystyle{amsplain}
\bibliography{arlist}

\providecommand{\bysame}{\leavevmode\hbox to3em{\hrulefill}\thinspace}
\providecommand{\MR}{\relax\ifhmode\unskip\space\fi MR }
\providecommand{\MRhref}[2]{%
  \href{http://www.ams.org/mathscinet-getitem?mr=#1}{#2}
}
\providecommand{\href}[2]{#2}
\begin{thebibliography}{1}

\bibitem{MR1801937}
Martin Aigner and G{\"u}nter~M. Ziegler, \emph{Proofs from {T}he {B}ook},
  second ed., Springer-Verlag, Berlin, 2001, Including illustrations by Karl H.
  Hofmann. \MR{MR1801937 (2001j:00001)}

\bibitem{MR2060634}
Va{\v{s}}ek Chv{\'a}tal, \emph{Sylvester-{G}allai theorem and metric
  betweenness}, Discrete Comput. Geom. \textbf{31} (2004), no.~2, 175--195.
  \MR{2060634 (2005c:52014)}

\bibitem{MR1194036}
J.~Csima and E.~T. Sawyer, \emph{There exist {$6n/13$} ordinary points},
  Discrete Comput. Geom. \textbf{9} (1993), no.~2, 187--202. \MR{1194036
  (94a:52015)}

\bibitem{MR0043485}
G.~A. Dirac, \emph{Collinearity properties of sets of points}, Quart. J. Math.,
  Oxford Ser. (2) \textbf{2} (1951), 221--227. \MR{0043485 (13,270c)}

\bibitem{MR0097014}
L.~M. Kelly and W.~O.~J. Moser, \emph{On the number of ordinary lines
  determined by {$n$} points}, Canad. J. Math. \textbf{1} (1958), 210--219.
  \MR{0097014 (20 \#3494)}

\end{thebibliography}

\end{document}